\documentclass[11pt,notitlepage]{article}
\usepackage{ams math} 
\usepackage{authblk} 
\usepackage{graphicx}
\usepackage{amssymb}
\usepackage{epstopdf}
\begin{document}

\title{Using Simulated Annealing to Factor Numbers}

\author[$^\dag$]{Eric Lewin Altschuler}
%
\author[$\dag$]{Timothy J. Williams}
\affil[$\dag$]{{\em Departments of Physical Medicine and Rehabilitation and Microbiology \& Molecular Genetics, Rutgers New Jersey Medical School}, {\tt altschel@njms.rutgers.edu}}
\affil[$\ddag$]{{\em Leadership Computing Facility, Argonne National Laboratory}, {\tt zippy@anl.gov}}



\maketitle

\begin{abstract}
Almost all public secure communication relies on the inability to factor large numbers. There is no known analytic or classical numeric method to rapidly factor large numbers. Shor\cite{shor1994algorithms} has shown that a quantum computer can factor numbers in polynomial time but there is no practical quantum computer that can yet do such computations. We show that a simulated annealing\cite{Kirkpatrick13051983} approach can be adapted to find factors of large numbers.
\end{abstract}

\section{Approach}\label{approach}

Simulated annealing is an approach inspired by statistical mechanics which by analogy views the values of a multivariable numeric problem as physical states of particles. Simulated annealing is very useful practically for problems such as the traveling salesperson problem and other so called NP-complete problems which have no known polynomial time numeric solutions. 

Let $N$ be the given number to be factored, and use $N$ in base 2 with $n$ digits. We seek numbers $A$ and $B$ (written in binary with respectively $a$ and $b$ digits) such that $A * B = N$ with $A$ the larger number (so $a>=b$). Now, $a + b = n$ or $n+1$. For a given $N$ there are at most $n-2$ possibilities for $\{a,b\}$. (Remember that $a>b$, and the leftmost digit of both A and B must be a 1.) Factor $A$, with $a$ binary digits, can have $1, 2, 3, \ldots, a$ 1s in its binary representation. Factor $B$, with $b$ binary digits, can have $1, 2, 3, \ldots, a$ 1s in its binary representation.

So, we formulate the factoring problem as follows: Given a binary number $N$ with $n$ digits, find binary numbers $A$ and $B$ with, respectively, $a$ and $b$ digits, of which $a^\prime$ and $b^\prime$ are 1, such that $A * B = N$.

The simulated annealing approach starts with a configuration of particles---here the binary digits of trial factors $A$ and $B$---and an energy $E$, here defined as:
\begin{equation}\label{energyDefinition}E = \sum_{i=1}^{n}\left\{ \begin{array}{ll}
					f(i) & \mbox{if $\{AB\}_i = N_i$} \\
					0 & \mbox{if $\{AB\}_i \ne N_i$}
				  \end{array}
			\right. \end{equation}
where $i$ indexes the digits in the binary representations and $f(i)$ is a function that increases monotonically with $i$. We have tested linear ($f(i)=i$) and quadratic ($f(i)=i^2$). This kind of function favors having as many of the binary digits of the product $AB$ match the corresponding binary digits of $N$ as possible, with weighting increasing for the higher bits. If all digits match exactly, $A$ and $B$ are factors of $N$, and the value of $E$ is maximal. If not, try a different configuration---here new possible factors $A^\prime$ and $B^\prime$---and look at the energy $E^\prime$ in this configuration. If the energy is increased then $A^\prime$ and $B^\prime$ are accepted as the new factors. Using the Metropolis algorithm,\cite{:/content/aip/journal/jcp/21/6/10.1063/1.1699114} occasionally $A^\prime$ and $B^\prime$ are accepted as the new trial factors even if $E^\prime<E$.  The chance of accepting the lower energy state is reduced as time goes on; by analogy, the problem cools according to an annealing schedule. From an initial ``temperature'' value $T_0$, given a cooling factor $F_c$, we iterate through $N_a$ annealing steps, each time reducing the temperature by the cooling factor. That is, moving from step $i$ to step $i+1$,
\begin{equation*}
T_{i + 1} = T_i * F_c
\end{equation*}

At each annealing step, try some number of configurations---rearrangements of the bits of binary representations of $A$ and $B$. We allow a number of different mechanisms to generate configurations. For one of the factors, choose from one of the following moves:

\begin{description}
\item[swap] Choose a pair of distinct bits (a 1 and a 0) at random, and swap them
\item[slide] Choose a random contiguous sub-sequence of the bits. Remove the rightmost bit, slide the remaining bits one to the right, then put the removed bit in the hole left behind.
\item[reverse] Choose a random contiguous sub-sequence of the bits and reverse its order.
\item[random] Randomly permute a random selection of the bits (generally a sparse selection, not a contiguous sequence).
\end{description}

For each configuration, test whether $A^\prime$ and $B^\prime$ multiply to $N$ (meaning they are factors of $N$). If so, repeat the whole annealing algorithm recursively on $A^\prime$ and $B^\prime$, finding successive sub-factors, until all endpoint numbers are are either prime factors (success) or the algorithm fails to factor one of them.

If $A^\prime B^\prime \neq N$, test whether whether the new energy of the system as defined in Eq. \ref{energyDefinition}, $E^\prime$, is greater than the energy of the previous configuration ($E$). If so, accept the new configuration as current and discard the previous one. If $E^\prime < E$, select a uniform random number between 0 and 1, $r$, and accept the new configuration if
\begin{equation*}
r < \exp\left(-\frac{E^\prime - E}{kT_i}\right)
\end{equation*}
where $T_i$ is the temperature at annealing step $i$ and $k$ is a constant analogous to the Boltzmann constant in physics.

This approach to factoring $N$ is nondeterministic, meaning that there is no guarantee of successfully finding all (or any) of the prime factors. However, the algorithm executes in polynomial time. If we know in advance that $N$ is {\em semiprime} (the product of two and only two prime factors), we stop the calculation once factors $A$ and $B$ are found such that $AB=N$.

The number of configurations tested per annealing step is a tunable parameter, as are the cooling factor, the number of annealing steps, and $k$. Note that our definition of $E$ implies that the optimum value is the largest $E$, not the smallest. The ``cooling'' actually leads to higher and higher values of $E$.

At this point, the reader may ask how we determined the number of bits in $A$ and $B$. If the number to be factored, $N$, has $n$ bits in its binary representation, its factors could have any number of bits from 2 to $n-1$. For a given $A$, having $a$ bits, we can compute the possible numbers of bits the other factor $B$ could have, given that $a + b = n$ or $n+1$. Our actual program explicitly factors out any prime factors up to (decimal) 1000; so if $N$ has $n$ bits, $A$ could have anywhere from 10 to $n-9$ bits. The corresponding $B$ could have either $b=n-a$ or $b=n-a+1$ bits.

Now the reader may ask how  we set the number of 1s in each binary number $A$ and $B$. We specify that the leftmost bit in both is 1 (that is, no leading 0 bits are allowed). With $a_1$ denoting the number of 1s in $A$, and $b_1$ denoting the number of 1s in $B$, we know that their range of values is
\begin{equation*}
\begin{array}{ll}
1\leq a_1 \leq a \\
1\leq b_1 \leq b
\end{array}
\end{equation*}
Our algorithm must try all $ab$ combinations. This scales roughly as $n^2$.

The whole algorithm then has a deep loop nest:
\begin{itemize}
\item{loop over the possible values of $a$}
\item{loop over the corresponding possible values of $b$}
\item{loop over the possible values of $a_1$}
\item{loop over the possible values of $b_1$}
\item{loop over all annealing steps (temperatures)}
\item{loop over configurations}
\end{itemize}

The scaling of the algorithm is then upper-bounded by
\begin{equation*}
n * (n-1) * n * (n-1) * N_a * N_c \approx n^4 * N_a * N_c
\end{equation*}
where $N_a$ is the number of annealing steps and $N_c$ is the number of configurations. Here we ignored the optimization of factoring out all small prime factors less than 1000, and approximated the maximum number if digits in $A$ and $B$ as $n$; their exact limits are lower, as detailed earlier. With $N_a$ and $N_c$ being constants, this is fourth order in $n$, the number of digits in the binary representation of the number we're factoring, $N$.

\section{Tests}

We tested our algorithm on numbers with up to 31 decimal digits (67 binary digits). Typically, the number of configurations per annealing step was set to
\begin{equation*}
M * \max(a,b)
\end{equation*}
where $a$ and $b$ are the number of binary digits in the current configuration of factors $A$ and $B$, and $M$ is an input parameter.

We constrain the acceptable configurations to only those for which the number of 1s in the binary representation of the product $AB$ is equal to the (known) number of 1s in the binary representation of $N$. Optionally, we constrain allowed ``bad'' moves to restrict the decrease in energy to be within a specified fraction of the current energy; this avoids ``really bad'' moves.  Optionally, we can retain the previous, higher energy configuration prior to an allowed ``bad'' move, run further configurations for a specified number of tries, then revert to the saved prior configuration if the energy has not evolved to exceed the prior energy. Table \ref{testCaseTable} and Table \ref{resultsTable} show results for a few test cases we successfully factored, and parameter values used.

\begin{table}
\centering
{
\small
\begin{tabular}{clll}
\hline
{\bf Case}  & {\bf $N$}                            & {\bf $n$}  & {\bf Factors} \\
\hline
A                & 99999989237606677        & 57            & $316227731*316227767$       \\
B                & 999999866000004473      & 60            & $999999929*999999937$       \\
C                & 9999999942014077477    & 64            & $3162277633*3162277669$   \\
D                & 99999980360000964323  & 67            & $9999999017*9999999019$   \\
\hline
\end{tabular}
\caption{{\small A few cases we tested. $N$ is the (semiprime) number to be factored. $n$ is the number of digits in the base-two representation of $N$.}}
\label{testCaseTable}
}
\end{table}

\begin{table}
\centering
{
\small
\begin{tabular}{ccrrrc}
\hline
{\bf Case}    & {\bf $N^*_a$}      & $F_c$ &$k$         & {\bf $N_c$}   & {\bf Time} \\
\hline
A                 & $41 \pm  62$       &0.997   & 63365   & 1,450,000     & $17 \pm 16$ \\
B                 & $13 \pm  8$         &0.997   & 73810   & 1,500,000     & $12 \pm 8$ \\
C                 & $144 \pm 161$    &0.997  & 89440    & 1,600,000     & $134 \pm 149$ \\
D                 & $89 \pm 40$        &0.997  &102510   & 1,700,000     & $96 \pm 45$ \\
\hline
\end{tabular}
\caption{{\small Results for a few cases we tested. $N^*_a$ is the number of annealing steps before finding the factors. $F_c$ is the temperature factor ($T$ is multiplied by this each annealing step. We always started with initial $T_0=F_c$.) $N_c$ is the maximum number of configurations tried each annealing step. Time is in runtime of the algorithm in minutes. In all cases, $N^*_a$ and Time are averaged over 5 runs, with the standard deviation indicated as error amounts. We used cost function $f(i)=i^2$ (see Eq. \ref{energyDefinition}).}}
\label{resultsTable}
}
\end{table}

The current implementation is in Python, using the {\tt bitarray} module for manipulating the binary representation. Python automatically handles arbitrary-precision integers correctly. We ran all the test cases on a laptop.

Since we knew the factors in advance for these test cases, we reduced the search space over numbers of bits and the number of 1s in the factors to be consistent with the known factors.  This is only an expedient, to be able to demonstrate the algorithm running on a single workstation. To run full tests for these and even larger numbers, with unknown factors, we are implementing a parallel version of the algorithm.

\section{Parallelism}

In the deep loop nest in Section \ref{approach}, there is ample parallelism to be exploited. All iterates of the outer four loops over numbers of digits and numbers of 1s in the factors can be computed independently, in parallel. This work scales as $n^4$, and is appropriate for distributed-memory parallelism with message passing. Memory requirements are very low, so replication of data is not a problem The parallel algorithm needs periodic, but infrequent, synchronization to test for success and proceed forward factorizing the factors, down to the final level of prime factors only. The loop over annealing steps is sequential, but the innermost loop over configurations can be parallelized; this is a candidate for thread parallelization and shared memory. Our initial target architecture will be an IBM Blue Gene/Q system, on which we have Python for the compute nodes.

\section{Conclusions}

We have shown the feasibility of using a nondeterministic optimization approach such as simulated annealing to work in a controlled and constrained manner toward prime factorization of large integers. Using a binary (base 2) representation allows for simple configuration-changing rules, and a simple energy cost function whose value is optimized. The algorithm has polynomial scaling in $n$, the numb of bits in the binary representation of the number to be factored. It is not guaranteed to find an exact solution, which is the only solution of interest in number factorization, but in practice we have found solutions for a wide range of numbers. It seems that the closer (further) the ratio of 1s to 0s in the binary representation of a factor of a semiprime number is to one, the harder (easier) it will be for a simulated annealing based method to factor the number. Thus, in picking semiprime numbers to use for encryption, a consideration should be the ratio of 1s to 0s in the two factors.

\bibliographystyle{unsrt}
\bibliography{SAFactor}

\begin{thebibliography}{1}

\bibitem{shor1994algorithms}
Peter~W Shor.
\newblock Algorithms for quantum computation: discrete logarithms and
  factoring.
\newblock In {\em Foundations of Computer Science, 1994 Proceedings., 35th
  Annual Symposium on}, pages 124--134. IEEE, 1994.

\bibitem{Kirkpatrick13051983}
S.~Kirkpatrick, C.~D. Gelatt, and M.~P. Vecchi.
\newblock Optimization by simulated annealing.
\newblock {\em Science}, 220(4598):671--680, 1983.

\bibitem{:/content/aip/journal/jcp/21/6/10.1063/1.1699114}
Nicholas Metropolis, Arianna~W. Rosenbluth, Marshall~N. Rosenbluth, Augusta~H.
  Teller, and Edward Teller.
\newblock Equation of state calculations by fast computing machines.
\newblock {\em The Journal of Chemical Physics}, 21(6):1087--1092, 1953.

\end{thebibliography}

\end{document}